\documentclass[11pt]{article}

\usepackage{amscd,amsmath, amssymb, fancyhdr, epsfig, color, tocloft, mathtools}
\usepackage{dutchcal}

\usepackage[backref=page]{hyperref}
\renewcommand*{\backrefalt}[4]{%
	\ifcase #1 (Not cited.)%
	\or        (Cited on page~#2.)%
	\else      (Cited on pages~#2.)%
	\fi}

\hypersetup{
	colorlinks   = true,
	citecolor    = magenta,
	linkcolor    = blue,
	urlcolor     = magenta	
}

\newcommand{\version}{version 1.0,\ \ May 27, 2022}

\renewcommand{\leftmark}{{\small Mall bundles and flat connections on Hopf manifolds}}

\makeatletter
\def\x@arrow{\DOTSB\Relbar}
\def\xlongequalsignfill@{\arrowfill@\x@arrow\Relbar\x@arrow}
\providecommand{\xlongequal}[2][]{%
	\ext@arrow 0099\xlongequalsignfill@{#1}{#2}}
\def\xlongrightarrowfill@{\arrowfill@\relbar\relbar\longrightarrow}
\newcommand{\xlongrightarrow}[2][]{%
	\ext@arrow 0099\xlongrightarrowfill@{#1}{#2}}
\makeatother

\numberwithin{equation}{section}

\newcommand{\dev}{{\operatorname{\sf dev}}}

\newcommand{\GL}{{\sf{GL}}}

\def\eqref#1{(\ref{#1})}

\newcommand{\goth}{\mathfrak}

\newcommand{\Z}{{\mathbb Z}}
\newcommand{\C}{{\mathbb C}}
\newcommand{\R}{{\mathbb R}}

\newcommand{\6}{\partial}
\def\1{\sqrt{-1}\,}

\newcommand{\restrict}[1]{{\left|_{{\phantom{|}\!\!}_{#1}}\right.}}

\newcommand{\cntrct}                
{\hspace{2pt}\raisebox{1pt}{\text{$\lrcorner$}}\hspace{2pt}}

\newcommand{\arrow}{{\:\longrightarrow\:}}

\newcommand{\calo}{{\cal O}}

\renewcommand{\bar}{\overline}
\renewcommand{\phi}{\varphi}
\renewcommand{\epsilon}{\varepsilon}
\renewcommand{\geq}{\geqslant}
\renewcommand{\leq}{\leqslant}

\newcommand{\im}{\operatorname{im}}
\newcommand{\End}{\operatorname{End}}

\newcommand{\Id}{\operatorname{Id}}

\renewcommand{\sup}{{\operatorname{{\sf sup}}}}

\newcommand{\const}{\operatorname{{\sf const}}}

\newcommand{\Aut}{\operatorname{Aut}}

\newcommand{\Sym}{\operatorname{\sf Sym}}

\renewcommand{\dim}{\operatorname{\sf dim}}

\newcommand{\rk}{\operatorname{rk}}

\newcommand{\Spec}{\operatorname{{\sf Spec}}}
\newcommand{\Alt}{\operatorname{Alt}}

\newcommand{\Pic}{\operatorname{Pic}}

\renewcommand{\Re}{\operatorname{Re}}

\newcounter{Mycounter}[section]
\newcounter{lemma}[section]
\renewcommand{\thelemma}{{Lemma \thesection.\arabic{lemma}}}
\newcommand{\lemma}{%
	\setcounter{lemma}{\value{Mycounter}}
	\refstepcounter{lemma}
	\stepcounter{Mycounter}
	{\noindent \bf \thelemma:\ }}

\newcounter{claim}[section]
\renewcommand{\theclaim}{{Claim \thesection.\arabic{claim}}}
\newcommand{\claim}{%
	\setcounter{claim}{\value{Mycounter}}
	\refstepcounter{claim}
	\stepcounter{Mycounter}
	{\noindent \bf \theclaim:\ }}

\newcounter{sublemma}[section]
\newcounter{corollary}[section]
\renewcommand{\thecorollary}{{Corollary \thesection.\arabic{corollary}}}
\newcommand{\corollary}{%
	\setcounter{corollary}{\value{Mycounter}}
	\refstepcounter{corollary}
	\stepcounter{Mycounter}
	{\noindent \bf \thecorollary:\ }}

\newcounter{theorem}[section]
\renewcommand{\thetheorem}{{Theorem \thesection.\arabic{theorem}}}
\newcommand{\theorem}{%
	\setcounter{theorem}{\value{Mycounter}}
	\refstepcounter{theorem}
	\stepcounter{Mycounter}
	{\noindent \bf \thetheorem:\ }}

\newcounter{conjecture}[section]
\newcounter{proposition}[section]
\renewcommand{\theproposition} {{Proposition \thesection.\arabic{proposition}}}
\newcommand{\proposition}{%
	\setcounter{proposition}{\value{Mycounter}}
	\refstepcounter{proposition}
	\stepcounter{Mycounter}
	{\noindent \bf \theproposition:\ }}

\newcounter{definition}[section]
\renewcommand{\thedefinition} {{Definition~\thesection.\arabic{definition}}}
\newcommand{\definition}{%
	\setcounter{definition}{\value{Mycounter}}
	\refstepcounter{definition}
	\stepcounter{Mycounter}
	{\noindent \bf \thedefinition:\ }}

\newcounter{example}[section]
\newcounter{remark}[section]
\renewcommand{\theremark}{{Remark \thesection.\arabic{remark}}}
\newcommand{\remark}{%
	\setcounter{remark}{\value{Mycounter}}
	\refstepcounter{remark}
	\stepcounter{Mycounter}
	{\noindent \bf \theremark:\ }}

\newcounter{problem}[section]
\newcounter{question}[section]
\makeatletter

\def\blacksquare{\hbox{\vrule width 5pt height 5pt depth 0pt}}
\def\endproof{\blacksquare}

\newcommand{\proof}{{\bf Proof: \ }}
\newcommand{\pstep}{{\bf Proof. Step 1: \ }}

\begin{document}

\begin{center}
{\Large\bf  Mall bundles and flat connections \\[3mm] on Hopf manifolds}\\[5mm]
{\large
Liviu Ornea\footnote{Liviu Ornea is  partially supported by Romanian Ministry of Education and Research, Program PN-III, Project number PN-III-P4-ID-PCE-2020-0025, Contract  30/04.02.2021},  
Misha Verbitsky\footnote{Misha Verbitsky is partially supported by
 the HSE University Basic Research Program, FAPERJ E-26/202.912/2018 
and CNPq - Process 313608/2017-2.\\[1mm]

}\\[4mm]

}

\end{center}

{\small
\hspace{0.15\linewidth}
\begin{minipage}[t]{0.7\linewidth}
{\bf Abstract} \\ 
A {\bf Mall bundle} on a Hopf manifold $H= \frac{\C^n \backslash 0}{\Z}$
is a holomorphic vector bundle whose pullback to $\C^n \backslash 0$ is trivial.
We define {\bf resonant} and {\bf non-resonant} Mall bundles, generalizing the notion of the
resonance in ODE, and prove that a non-resonant Mall bundle always admits
a flat holomorphic connection. We use this observation to prove a version of 
Poincar\'e-Dulac linearization theorem, showing that any non-resonant
invertible holomorphic contraction of $\C^n$ is linear in appropriate
holomorphic coordinates. We define
the notion of {\bf resonance} in Hopf manifolds, and show that all non-resonant
Hopf manifolds are linear; previously, this result was obtained by Kodaira using
the Poincar\'e-Dulac theorem.
\end{minipage}
}

\tableofcontents

\section{Introduction}

\subsection{Mall bundles: history and definition}

Let $f:\; M \arrow M$ be a continuous map fixing $x\in M$.
We say that $f$ is {\bf a contraction centered in $x$}
if for any compact set $K \subset M$ and any
open set $U\subset M$ containing $x$,
a sufficiently high power of $f$ takes $K$ to $U$.

Hopf manifolds are quotients of $\C^n \backslash 0$ by a
holomorphic contraction centered at $0$.
In his seminal work \cite{_Mall:Contractions_}, Daniel Mall 
computed the cohomology
of a holomorphic vector bundle on a Hopf manifold
such that its pullback to $\C^n \backslash 0$ can be
extended to a holomorphic vector bundle on $\C^n$.
We call such vector bundles {\bf the Mall bundles}
(\ref{_Mall_bundle_Definition_}). The same argument
was earlier applied by A. Haefliger, \cite{_Haefliger:Hopf_}.

By the Oka-Grauert homotopy principle
(\cite[Theorem 5.3.1]{_Forstneric:Oka_book_}),
any vector bundle on $\C^n$ is trivial.
Therefore, one could define the Mall bundles
as holomorphic vector bundles on $H$ such that their
pullback to $\C^n \backslash 0$ is trivial,
\ref{_Mall_via_trivial_Remark_}.

Since 1990-ies, Mall's theorem was explored and
generalized in many different directions; see,
for example \cite{_Libgober_} and \cite{_Gan_Zhou_}.

\subsection{Mall theorem and its applications}

The main utility of Mall bundles is
Mall theorem, which allows one to compute
the cohomology of a Mall bundle in terms of 
holomorphic sections. This theorem says
that the cohomology $H^i(H, B)$ of a Mall bundle
$B$ on a Hopf manifold $H=\frac{\C^n\backslash 0}{\Z}$
vanish for $i\neq 0, 1, n-1, n$, and that
$\dim H^0(H, B)= \dim H^1(H, B)$ and
$\dim H^{n-1}(H, B)= \dim H^n(H, B)$.
Serre duality gives an isomorphism
$H^n(H, B)= H^0(H, B\otimes K_H)^*$,
hence all cohomology of $B$ can
be expressed in terms of the holomorphic sections.

We use this observation on cotangent bundle
and the tensor bundles of form $(\Omega^1H)^{\otimes n} \otimes_{\calo_H} TH$.
It is not hard to see (\ref{_Mall_flat_connection_Theorem_}, Step 2)
that the obstruction to existence
of a holomorphic connection belongs to $H^1(\Omega^1 H\otimes \End(B))$.
Therefore, a Mall bundle admits a holomorphic connection whenever
$\dim H^0(\Omega^1 H\otimes \End(B))=0$. A Mall bundle
which satisfies this condition is called {\bf non-resonant}.
A non-resonant Mall bundle has a unique holomorphic
connection, which is a posteriori flat (\ref{_Mall_flat_connection_Theorem_}).

This is used to construct a torsion free flat affine
connection on a non-resonant Hopf manifold, proving the
Poincar\'e linearization theorem.

It would be interesting to access the existence 
of holomorphic connections on a resonant Mall bundle.
We could not find a Mall bundle which does not
admit a holomorphic connection. It is clear that
a non-linearizable Hopf manifold cannot admit
a flat, torsion-free holomorphic connection 
(\ref{_flat_Hopf_is_linear_Theorem_}).
However, it is not clear whether the connection
with non-zero holomorphic curvature can exist.

\subsection{Mall bundles and holomorphic connections}

In this paper we explore the geometric properties of the
Mall bundles.

Let $\gamma$ be an invertible holomorphic contraction of
$\C^n$ centered in 0, and $H= \frac{\C^n \backslash  0}{\langle \gamma\rangle}$
the corresponding Hopf manifold. 
Since the pullback of a Mall bundle to $\C^n \backslash 0$
is a trivial bundle, the category of Mall bundles
is equivalent to the category of $\gamma$-equivariant
bundles on $\C^n$.

Let $B$ be a complex vector bundle on a complex manifold $M$,
and $\nabla:\; B \arrow B \otimes \Lambda^1 M$
 a flat connection. The Hodge component $\bar\6:=\nabla^{0,1}$
is a holomorphic structure operator on $B$. By Koszul-Malgrange theorem (\cite[Chapter I, Proposition
3.7.]{_Kobayashi_Bundles_}, \cite{_Koszul_Malgrange_}), 
the sheaf ${\cal B}:=\ker\bar\6$ is a holomorphic vector bundle on
$M$, with ${\cal B} \otimes_{\calo_M} C^\infty M = B$.
In this situation we say that the flat
connection $\nabla$ {\bf is compatible with the
holomorphic structure on $B$}.

There are many examples of Mall bundles arising from the
geometry of Hopf manifolds. All tensor bundles, all line
bundles, and all extensions of Mall bundles are also Mall
(\ref{_Mall_exa_Proposition_}).
In many of those examples, the equivariant action of $\gamma$
on $B$ preserves a flat connection on $B$ 
(for the line bundles, it follows from 
\ref{_Picard_Hopf_Proposition_}). In other words,
these Mall bundles are obtained from flat bundles
by taking the $(0,1)$-part of the connection.
It turns out that this situation is quite general, and 
an arbitrary Mall bundle admits a compatible flat
connection when the so-called ``non-resonance'' condition
is satisfied.\index{resonance}

This can be explained as follows. Recall that 
{\bf a holomorphic connection} (\ref{_holo_conne_Definition_}) on a holomorphic vector\index{connection!holomorphic}
bundle $B$ is a holomorphic differential operator
$\nabla:\; B \arrow B \otimes \Omega^1 M$ satisfying the
Leibniz rule, $\nabla(fb) = f \nabla(b) + df \otimes \nabla(b)$
for any holomorphic $f\in \calo_U$ and any $b \in H^0(U, B)$. 
We want to construct a holomorphic connection on a Mall
bundle on a Hopf manifold; this is equivalent to having
a $\gamma$-equivariant connection on its pullback to
$\C^n\backslash 0$ considered as a $\gamma$-equivariant
vector bundle.

Consider
the space ${\cal A}$ of all holomorphic connections on a trivial
$\gamma$-equivariant holomorphic vector
bundle $R$ on $\C^n$. This is an affine  space
modeled on the vector space $H^0(\Omega^1 \C^n\otimes_{\calo_{\C^n}} \End(R))$,
and the equivariant action defines an affine endomorphism
of ${\cal A}$. The linearization $\rho$ of this action is a compact
endomorphism of $H^0(\Omega^1 \C^n\otimes_{\calo_{\C^n}} \End(R))$
considered as a topological vector space with $C^0$-topology
(\ref{_eigenvalues_diff_forms_contraction_Lemma_}).
If $\rho$ has all eigenvalues with absolute value $< 1$,
Banach fixed point theorem would imply that $R$
admits a $\gamma$-equivariant holomorphic connection.
In fact, it would suffice to check that
all eigenvalues $\lambda_i$ of $\rho$ 
are not equal to 1.

The ``resonance'' is a property of the
eigenvalues of the $\gamma$-equivariant action;
a $\gamma$-equivariant vector bundle $R$
on $\C^n$ {\bf has resonance} when the eigenvalues
of $D\gamma$ on the fiber $R \restrict 0$
are $\beta_1, ..., \beta_m$, the eigenvalues
of $D\gamma$ on $T_0 \C^n$ are $\alpha_1, ..., \alpha_n$,
and there exists a relation of the form\index{resonance!of a bundle}
$\beta_p = \beta_q \prod_{i=1}^n \alpha_i^{k_i}$, with
all $k_i$ non-negative integers, and $\sum_i k_i > 0$,
for some $\beta_p, \beta_q$, which are not necessarily distinct.

Let $B$ be a vector bundle on a Hopf manifold, and
$R$ the extension of its pullback to $\C^n$. We prove that
$R$ has no resonance if and only if $H^0(H,
\Omega^1_H\otimes \End(B))=0$ (\ref{_B_on_Hopf_resonant_via_sections_Corollary_}).
We also prove that any non-resonant holomorphic vector bundle
on a Hopf manifold admits a flat connection compatible
with the holomorphic structure (\ref{_Mall_flat_connection_Theorem_}).

\subsection{Flat affine structures and the development map}

The notion of resonance is classical, and harks back to Poincar\'e,
Latt\`es and Dulac, who discovered the resonance while working on
the normal forms of ordinary differential equations.
In the modern language, they were looking at the normal form
of a real analytic or complex analytic vector field which has a simple
zero at a given point.

The ``normal form'' is a classical notion, which is roughly equivalent,
in the modern language, to ``the moduli space'', but includes a more 
explicit description in terms of coordinates. 

For vector fields without zeroes, the normal form is very simple:
in an appropriate coordinate system, this vector field takes
the form $\frac d{dx_1}$; this is called ``straightening of a vector field''.
This result follows directly from the Peano and Picard theorems on
the existence of solutions of ODE.

The normal form theorem for a vector field with a simple zero 
is known as the {\em Poincar\'e-Dulac theorem}, \cite{_Lattes_,_Dulac_,_Arnold:ODE+_}. 
We give a general outline of this theory, following \cite{_Enc_Math:Poincare_Dulac_}. 

Let $\dot x = A(x) + u(x)$ be a formal
differential equation, where $x(t)\in \C^n$ is a time-dependent
point in $\C^n$, $A$ a non-degenerate linear operator, and $u(x)$ 
a Taylor series starting from the second order terms. It is said that {\bf $A$ has a resonance} if
there is a relation of the form $\lambda_i = \sum_{j=0}^n m_j\lambda_j$,
where $m_j \in \Z^{\geq 0}$ and $\sum_{j=0}^n m_j\geq 2$.\footnote{Later in this
chapter, we redefine this notion in such a way that this additive
relation becomes multiplicative, $\lambda_i = \prod_{j=0}^n \lambda_j^{m_j}$;
this is done because we work with holomorphic contractions and not with the vector fields.}

If $A$ has no resonance, then the normal form of 
this vector field is very simple: in appropriate coordinates
$y_1, ..., y_n$ it can be written as $\dot y = A(y)$.
If $A$ has a resonance, the vector field has a normal form, 
which is written in a coordinate system  $y=(y_1, ..., y_n)$ as follows. 
Choose $y_i$ in such a way that $A$ is upper triangular
in the basis $\frac d{dy_i}$, and the diagonal
terms corresponding to $\frac d{dy_i}$ are $\lambda_i$.
Then the equation $\dot x = A(x) + w(x)$
has normal form $\dot y = A(y)+ \sum e_i w_i(y)$, where $e_i= \frac d{dy_i}$ is the
coordinate vector field, and $\sum e_i w_i(y)$ is a Taylor series
obtained as a sum of {\em resonant monomials}. A coordinate
monomial $e_i \prod_{j=1}^n y_j^{m_j}$ is {\bf resonant}
if $\lambda_i = \sum m_j \lambda_j$.

In general, it is hard to achieve convergence for these formal sums,
even when the differential equation is analytic. However, if $e^A$
is a contraction, the convergence is automatic, because the number
of resonant monomials is finite. Indeed, $e^A$ is a contraction
if and only if $\Re \lambda_i <0$ for all $i$, and 
the equation $\Re \lambda_i = \sum m_j \Re \lambda_j$, $m_i \in \Z^{\geq 0}$,
implies that $m_i \leq \max_{j,l}\frac {\Re \lambda_j}{\Re \lambda_l}$.

A similar result is true for germs of 
biholomorphic contractions, due to S. Sternberg (\cite{_Sternberg_contraction_}).
However, in this case one should replace the linear resonance
by multiplicative, $\lambda_i = \prod_{j=0}^n \lambda_j^{m_j}$,
as in \ref{_resonant_matrix_Definition_}.

An invertible holomorphic contraction gives rise to a Hopf manifold,
and Sternberg's theorem can be interpreted as a structure theorem
about Hopf manifolds; this is how Kodaira used it in \cite{_Kodaira_Structure_II_}.\index{theorem!Sternberg}

In this paper, we use the flat connection inherent on\index{vector bundle!Mall}
Mall bundles to give a new proof of the non-resonant part of the
Poincar\'e-Dulac theorem. Let $\gamma$ be a germ of an invertible 
biholomorphic contraction of $\C^n$ with center in 0. We say that
$\gamma$ {\bf is non-resonant} if the differential $D_0 \gamma\in \End(\C^n)$
is a non-resonant matrix. Let $H$ be the corresponding Hopf manifold;
then the tangent bundle $TH$ is non-resonant, which is equivalent
to $H^0(H, \Omega^1H \otimes \End(TH))=0$ (\ref{_B_on_Hopf_resonant_via_sections_Corollary_}).
This immediately implies that the flat connection in $TH$, given by
\ref{_Mall_flat_connection_Theorem_}, is torsion-free.

To prove the Poincar\'e-Dulac linearization theorem, we need to
find the coordinates on $\C^n$ in which $\gamma$ is linear.
To produce the flat coordinates, we use the developing map
defined in the framework of flat affine geometry
(or, more generally, in Cartan geometries).

A manifold $M$ is called {\bf affine}, or {\bf flat affine},
if it is equipped with an atlas of open sets, identified with
open subsets in $\R^n$, with the transition functions
affine. This is equivalent to having a torsion-free
flat affine connection on $M$
(\ref{_flat_affine_via_connection_Proposition_}).
The study of compact flat affine manifolds is 
ongoing, with many conjectures still open.
We refer to \cite{_Abels:survey_} for more details
and open questions.

%

Each flat affine manifold
is equipped with the natural flat, torsion-free connection $\nabla$.
Using this connection, the developing map
can be defined as follows. Assume that $(M, \nabla)$
is a simply connected flat affine manifold.
Let $\theta_1, ..., \theta_n$ be parallel 1-forms
which trivialize the bundle $T^*M$. Since $\nabla$ is
torsion-free, all the forms $\theta_i$ are closed;
however, $H^1(M)=0$, which implies that $\theta_i$ are
exact, $\theta_i= dz_i$. The map $m \arrow (z_1(m), ..., z_n(m))$
is called {\bf the developing map}.

We can also understand the developing map
in terms of the geodesics, as the inverse of the {\em exponential map}.
The exponential map in this context is a map
taking a tangent vector to the point at time 1 on the geodesic 
tangent to this vector in time 0. This definition is equivalent to the one
given above (\ref{_dev_for_complete_affine_Theorem_}).
A flat affine structure is {\bf complete}
when the geodesic equation $\nabla_{\dot \gamma_t} \dot \gamma_t$ can be solved
for all $t\in \R$ and all initial conditions $\gamma_0\in M$, $\dot \gamma_0\in T_{\gamma_0}M$.

The completeness condition is tricky and counter-intuitive;
indeed, even a compact flat affine manifold is not necessarily complete.
A textbook example of a non-complete flat affine manifold is a real linear Hopf
manifold $H$, obtained as a quotient of $\R^n \backslash 0$
by a linear contraction. This manifold is compact, but 
its universal cover is $\R^n \backslash 0$, and the 
developing map is an open embedding
$\R^n \backslash 0\hookrightarrow \R^n$. To obtain a non-complete geodesic,
one needs to start from a geodesic in $\R^n$ passing through
0; its image in $H$ is manifestly non-complete.

Let $M$ be a complete, simply connected affine manifold.
Then the developing map $\dev:\; M \arrow \R^n$
is an isomorphism of affine manifolds. This is a classical
result by Auslander-Markus (\cite{_Auslander_Markus:holonomy_}) that we 
prove in \ref{_dev_for_complete_affine_Theorem_}.

For our present purposes, we need a variation
of this result, which ultimately implies the non-resonant
case of the Poincar\'e-Dulac theorem. From \ref{_Mall_flat_connection_Theorem_},
it follows that any non-resonant 
Hopf manifold $M= \frac{\C^n \backslash 0}{\langle A \rangle}$ is equipped with a unique
torsion-free flat affine connection compatible with the
complex structure. However, it is not complete, as we explained above.
We prove that this flat affine connection lifted to the universal 
covering $\C^n \backslash 0$ of $M$ can be extended to 0,
resulting in a complete flat affine structure on $\C^n$.
The corresponding developing map puts flat affine coordinates
on $\C^n$, and the contraction $A$ is affine in these coordinates, 
hence linear. This gives a new proof of the non-resonant
case of the Poincar\'e-Dulac theorem.

\section{Preliminaries on Banach spaces}


We gather here what is needed in the sequel concerning
compact operators on Banach spaces and the Riesz-Schauder
theorem.

\hfill

Recall that 
a  subset $X$ of a topological space $Y$ is called {\bf 
	precompact},
or {\bf relatively compact in $Y$}, 
if its closure is compact. 

\hfill

\definition \label{_bounded_set_Definition_}
A subset
$K\subset V$ of a topological vector space
is called {\bf bounded} if for any 
open set $U\ni 0$, there exists a number $\lambda_U\in
\R^{>0}$ such that $\lambda_U K \subset U$.

\hfill

\definition \label{_compact_operator_Definition_}
Let $V, W$ be topological vector spaces, A continuous operator $\phi:\; V
\arrow W$ is called {\bf  compact} if the image of
any bounded set is precompact.

\hfill

Now Montel's theorem (\cite[Lemma 1.4]{_Wu:Montel_}) can be restated as follows:

\hfill

\claim \label{idcomp}
Let $V=H^0(\calo_M)$ be the space of holomorphic functions on a complex
manifold $M$ with $C^0$-topology. Then any 
bounded subset of $V$ is precompact. In this case, the
identity map is a compact operator. \endproof

\hfill

\remark
It is not hard to deduce from Montel theorem 
that the space of bounded holomorphic functions on $M$
is Banach (that is, complete as a metric space) 
with respect to the $\sup$-norm.

The following theorem can be used to obtain a version of the Jordan normal
form for a compact operator on a Banach space. Recall
that the {\bf spectrum} of a linear operator $F$ is the set of all
$\mu\in \C$ such that $F - \mu\Id$ is not invertible.

\hfill

\theorem{ (Riesz-Schauder, \cite[Section 5.2]{friedman})}\\
\label{_Riesz_Schauder_main_Theorem_}
Let $F:\; V \arrow V$ be a compact operator on a Banach
space.  Then the spectrum $\Spec F\subset \C$ is 
compact and discrete outside of $0\in \C$. Moreover, 
for each non-zero $\mu \in \Spec F$, there exists a sufficiently
big number $N\in \Z^{\geq 0}$ such that for each $n>N$ {one has 
	\[ V= \ker(F-\mu\Id)^n \oplus \overline{\im (F-\mu\Id)^n},
	\]
	where $\overline{\im (F-\mu\Id)^n}$ is the closure of the image.}
Finally, the space $\ker(F-\mu\Id)^n$ is finite-dimensional. \endproof

\hfill

\remark\label{_root_space_RS_Remark_}
Recall that   {\bf the root space of an operator 
	$F\in \End(V)$, associated with an eigenvalue $\mu$}, 
is $\bigcup_{n\in \Z} \ker(F-\mu\Id)^n$. A vector $v$ is called {\bf a root vector} for the operator $F$ if $v$ lies in 
a root space of $F$, for some eigenvalue $\mu\in \C$. 
In the finite-dimensional case, the root spaces
coincide with the  Jordan cells of the corresponding 
matrix. Then 
\ref{_Riesz_Schauder_main_Theorem_} can be reformulated
by saying that any compact operator $F\in \End(V)$ admits a 
Jordan cell decomposition, with $V$ identified with
a completed direct sum of the root spaces, which are
all finite-dimensional; moreover, the eigenvalues
of $F$ converge to zero.

\hfill

In the sequel, we shall use the following corollary of the
Riesz-Schauder theorem. 


\hfill

\theorem\label{rs}
Let  $F:\; V \arrow V$ be a compact operator on a Banach space.
Then the space generated by the root vectors is dense in $V$.
\endproof



%



%

%

\section{Dolbeault cohomology of Hopf manifolds}
\subsection{Computation of $H^{0,p}(H)$ for a Hopf manifold}

Dolbeault cohomology of Hopf manifolds is a classical subject,
but we could not find the computation for the general case.
For the classical Hopf manifold (a quotient of $\C^n \backslash 0$
by a constant times identity), an answer is given in 
\cite{_Ise_}. More general Hopf surfaces were defined and classified by Kodaira
(\cite{_Kodaira_Structure_II_}); he 
computed some of their cohomology in \cite{_Kodaira_Structure_III_}.
For a diagonal Hopf manifold, the Dolbeault cohomology
was computed by D. Mall
(\cite{_Mall_}). For a reference to other special cases of this
theorem, see \cite{_Ramani_Sankaran_,_MO:Hopf_cohomology_}.

\hfill

We recall the statement of Mall's result:

\hfill

\theorem\label{_Mall_cohomology_Theorem_} (\cite{_Mall:Contractions_})\\
Let $\pi:\;  \C^n\backslash 0\arrow H$,  $n \geq 3$, be the
universal cover of a Hopf manifold,  
$j:\; \C^n\backslash 0\hookrightarrow \C^n$
the standard embedding map, and 
$B$ a holomorphic vector bundle over $H$
such that $j_* \pi^* B$ is a locally trivial 
coherent sheaf on $\C^n$.%
\footnote{By Oka-Grauert 
homotopy principle, \cite[Theorem 5.3.1]{_Forstneric:Oka_book_},
any holomorphic vector bundle on $\C^n$ is trivial;
thus, instead of local triviality, we could assume 
that $j_* \pi^* B$ is a trivial vector bundle.}
Then $\dim H^0(H, B)=\dim H^1(H,B)$, and this group
is equal to the space of $\Z$-invariant sections
of $j_* \pi^* B$. Moreover, $H^i(H, B) =0$ for all
$i$ such that $1 < i < n-1$. \endproof

\hfill

Using Mall's theorem, one can easily compute the $(0,*)$-part of the 
Dolbeault cohomology of a Hopf $n$-manifold $H$. 
Indeed, $H^{0,1}(H)= H^1(\calo_H)$ has 
the same rank as $H^0(\calo_H)=\C$, and $H^{0,i}(H)=0$
for $1 < i < n-1$. By Serre duality, $H^{0,n-i}(H)= H^i(K_H)^*$,
where $K$ is the canonical bundle, hence to prove that 
$H^{0,n-1}(H)= H^1(K_H)^*$ and $H^{0,n}(H)= H^0(K_H)^*$
vanishes, it would suffice to show that the canonical
bundle on a Hopf manifold has no holomorphic sections.

\hfill

\theorem\label{_Dolbeault_for_Hopf_Theorem_}
Let $H= (\C^n \backslash 0)/\Z$ be a Hopf manifold,
that is, a quotient of $\C^n \backslash 0$ by a holomorphic contraction,
Then $H^i(\calo_H)=0$ unless $i = 0,1$, and
$\rk H^1(\calo_H)=\rk H^0(\calo_H)=1$.

\hfill

\proof
By Mall's theorem 
(\ref{_Mall_cohomology_Theorem_}), $\rk H^1(H,\calo_H)=\rk H^0(H,\calo_H)$,
and $H^i(H,\calo_H)=0$ for $1 <i < n-1$.
However, $\rk H^0(H,\calo_H)=1$ because $\calo_H$ is a trivial line bundle.
To finish the proof, it remains only to show that
$H^{n-1}(H,\calo_H)$ and $H^{n}(H,\calo_H)$ vanish.
By Serre duality, these two spaces are dual to $H^0(H, K_H)$ and
$H^1(H, K_H)$, which have the same rank by Mall's theorem again. 
It remains only to prove that
$H^0(H, K_H)=0$ for any Hopf manifold. 

Suppose that $\eta$ is a non-zero element in $H^0(H, K_H)=0$;
we consider $\eta$ as a holomorphic volume form. Then
$\mu:= \eta\wedge \bar \eta$ is a measure on $H$, which is
strictly positive outside of the zero divisor of $\eta$.

Consider the measure $\pi^* \mu:=\pi^* \eta\wedge \pi^* \bar\eta$ 
on $\C^n \backslash 0$. Since $\pi^* \eta$ is $\Z$-invariant, 
the measure $\pi^* \mu$ is $\Z$-invariant as well.
The canonical bundle of $\C^n$ is trivial, hence,
by Hartogs theorem, $\pi^*\eta$ can be extended to
a holomorphic section $j_* \pi^*\eta$ of $K_{\C^n}$.
Denote by $j_* \pi^* \mu$ the measure
$j_*\pi^* \mu:=j_*\pi^* \eta\wedge j_*\pi^* \bar\eta$.
This measure is finite on compacts, and is preserved
by the contraction $\gamma:\; \C^n \arrow \C^n$. 
This is impossible, unless $j_*\pi^* \mu=0$, because any bounded set
is mapped inside a given compact neighbourhood
of 0 by a sufficiently big power of $\gamma$.
This implies that $\eta=0$, hence $H^0(H, K_H)=0$.
\endproof

\subsection{Holomorphic differential forms on Hopf manifolds}

The results of this subsection generalize the vanishing 
$H^0(H, K_H)=0$ given in the proof of \ref{_Dolbeault_for_Hopf_Theorem_}.
It turns out that all holomorphic differential forms on Hopf
manifolds vanish. 

\hfill

In \cite{_Ise_} the vanishing of differential forms was proven
for the classical Hopf manifold $\frac{\C^n \backslash 0}{\lambda \Id}$;
we could not find other results in the literature, though the question
seems to be elementary and classical.


%
%
%

\hfill

We start from the following lemma.

\hfill

\lemma \label{_eigenvalues_diff_forms_contraction_Lemma_}
Let $\gamma:\; \C^n \arrow \C^n$ be an invertible holomorphic contraction
centered in 0, and $D \subset \C^n$  an open set such that
$\gamma(D)$ is precompact in $D$. Choose an Hermitian
metric on $\C^n$, and define the norm on the space
$H^0_b(D, \Omega^1D)$ of bounded holomorphic $1$-forms 
as $\|\eta\|:= \sup_{x\in D} |\eta_x|$.\footnote{%
By Montel's theorem, $H^0_b(D, \Omega^1D)$
with this norm is a Banach space.}  
Then the operator 
$\gamma^*:\; H^0_b(D, \Omega^1D)\arrow  H^0_b(D, \Omega^1D)$
is compact, and all its eigenvalues are smaller than 1 in absolute value.

\hfill

For the proof, we shall need the following result about the compactness of contraction operators:

\hfill

\theorem \label{_contra_compact_Theorem_}\cite[Theorem 2.14]{ov_indam}\\
Let $X$ be a complex variety, and 
$\gamma:\; X \arrow X$ a holomorphic contraction
centered in $x\in X$ such that $\gamma(X)$ is precompact. 
Consider the Banach space $V=H^0_b(\calo_X)$ of bounded holomorphic
functions with the sup-norm, and let $V_x\subset V$
be the space of all $v\in V$ vanishing in $x$. Then the operator $\gamma^*:\; V \arrow V$
is compact, and the eigenvalues  
of its restriction to $V_x$ are strictly smaller than 1
in absolute value.\footnote{Since
$\gamma^*$ maps constants to constants identically, one cannot expect
that the eigenvalues of $\gamma^*$
satisfy $|a_i|< 1$ on $V$. However, if we add a condition
which excludes constants, such as $v(x)=0$, we immediately
obtain $|a_i|< 1$.} \endproof

\hfill

%

{\bf Proof of \ref{_eigenvalues_diff_forms_contraction_Lemma_}:}
The operator 
$\gamma^*:\; H^0_b(D, \Omega^1D)\arrow  H^0_b(D, \Omega^1D)$
is compact by an argument which follows from Montel's theorem 
(\ref{_contra_compact_Theorem_}). 

The differential of $\gamma$
acts with all eigenvalues $|a_i|<1$ on $T_0\C^n$,
by \ref{_contra_compact_Theorem_}. 
Then $\gamma^*$ (the pullback operator on differential forms)
acts on $T^*_0\C^n$ with all eigenvalues $|a_i|<1$.
Consider the Taylor expansion of a function $f$ in 0.
The chain rule and the estimate of the eigenvalues
of $\gamma^*$ on $T^*_0\C^n$ imply that $\gamma^*$ acts
on the non-constant Taylor coefficients of $f$ with all eigenvalues
$|a_i|<1$.

Let $a$ be an eigenvalue of the
compact operator 
$$\gamma^*:\; H^0_b(\Omega^1D) \arrow H^0_b(\Omega^1D).$$
It remains to show that $|a|<1$.

The eigenvalues of $\gamma^*$ on the vector space
$\Omega^1\C^n\restrict 0$ are powers of its action on
the cotangent bundle, which are all smaller than 1 in absolute value,
because $\gamma$ is a contraction. This implies that the eigenvalues
of $\gamma^*$ on the degree 0 term of the Taylor
expansion of $\eta$ are also smaller than 1 in absolute value.
Summarizing the above estimates, we obtain that
$\lim_n (\gamma^*)^n \eta=0$, hence $\gamma^*$
acts on differential forms with all eigenvalues
smaller than 1 in absolute value.
\endproof

\hfill

\proposition\label{_holo_forms_Hopf_vanish_Proposition_}
Let $H$ be a Hopf manifold, and $B$ a tensor power of $\Omega^1H$,
$B= (\Omega^1H)^{\otimes l}$. Then all holomorphic sections of $B$ vanish.

\hfill

\pstep
The universal cover of the Hopf manifold is $\C^n \backslash 0$
equipped with the free and properly discontinuous holomorphic $\Z$-action.
Choose precompact fundamental domains  $M_i$ of the $\Z$-action,
with the generator of $\Z$ mapping $M_i$ to $M_{i+1}$.
Then we may assume that
$\tilde M_0 :=\{0\}\cup \bigcup_{i\leq 0} M_i$
is a bounded neighbourhood of 0 in $\C^n$.

Let $\gamma\in \Aut(\C^n)$ be the generator of the $\Z$-action
acting on $\C^n$ as a contraction.
By  \ref{_eigenvalues_diff_forms_contraction_Lemma_}, 
the action of $\gamma^*$ on $B$ defines a compact operator 
\[ H^0_b(\tilde M_0, B) \arrow H^0_b(\tilde M_0, B),\]
with all eigenvalues smaller than 1.

Given a holomorphic differential form $\alpha$ on $\C^n$,
we restrict it to $\tilde M_0$ and observe that it is bounded
because $\tilde M_0$ is precompact.
Using the compactness of $\gamma^*$-action and the estimate of its eigenvalues
(\ref{_eigenvalues_diff_forms_contraction_Lemma_}), we show that
the norm of  $(\gamma^n)^*\alpha$ on $\tilde M_0$ converges to 0.

\hfill

{\bf Step 2:} 
Let $\pi:\; \C^n\backslash 0 \arrow H$ be the universal cover of a Hopf
manifold, and $\beta\in H^0(H, B)$ a section of $B$.
Any holomorphic differential form on $\C^n \backslash 0$
can be extended to $\C^n$, when $n>1$, by Hartogs theorem.
Therefore, the section 
$\pi^*\beta$ can be holomorphically extended to 
0, and this extension is
$\gamma^*$-invariant. By the argument in Step 1,
$\lim_n (\gamma^n)^*\alpha=0$
for any holomorphic form $\alpha$; applying
this to $\pi^*\beta$, we conclude that $\beta=0$.
\endproof

\section{Mall bundles on Hopf manifolds}

\subsection{Mall bundles: definition and examples}

We define {\bf Mall bundles} on a Hopf manifold
as bundles which satisfy the assumptions of 
\ref{_Mall_cohomology_Theorem_}.

\hfill

\definition\label{_Mall_bundle_Definition_}
Let $\pi:\;  \C^n\backslash 0\arrow H$ be the
universal cover of a Hopf manifold,  
$j:\; \C^n\backslash 0\hookrightarrow \C^n$
the standard embedding map, and 
$B$ a holomorphic vector bundle over $H$
such that $j_* \pi^* B$ is a locally trivial
coherent sheaf on $\C^n$, that is, a holomorphic vector bundle. 
Then $B$ is called {\bf a Mall bundle}.

\hfill

Note that any holomorphic vector bundle on $\C^n$
is trivial, as follows from the Oka-Grauert homotopy principle
(\cite[Theorem 5.3.1]{_Forstneric:Oka_book_}).\index{theorem!Oka-Grauert homotopy principle}

\hfill

\remark\label{_Mall_via_trivial_Remark_}
Let $B$ be a Mall bundle on a Hopf manifold $H$.
Then its pullback  $\pi^* B$ to $\C^n \backslash 0$ is extended
to a trivial holomorphic bundle $\hat B$, hence 
$\pi^* B=\hat B\restrict{\C^n \backslash 0}$ is trivial.
Conversely, if  $\pi^* B$ is trivial on $\C^n \backslash 0$, it 
can be extended to a bundle on $\C^n$, hence it is Mall. One could define
Mall bundles as holomorphic bundles on $H$ such that
$\pi^* B$ is a trivial bundle on $\C^n \backslash 0$.

\hfill

Since $j_*$ commutes with direct sums
and tensor products, the following observation is clear.

\hfill

\claim\label{_tensor_direct_component_Mall_Claim_}
A tensor product of Mall bundles, a direct sum
of Mall bundles, and any direct sum component
of a Mall bundle is again a Mall bundle.

\hfill

\proof For direct sums and tensor products the statement
is obvious, because these operations commute with the
functor $j_* \pi^*$. For a direct sum component, 
consider a coherent sheaf ${\cal F}$ on a complex manifold $M$, and let
${\cal F}|_x:= {\cal F}\otimes_{\calo_M} (\calo_M/{\goth m}_x)$.
This is a finite-dimensional space; if ${\cal F}$ is a vector bundle,
${\cal F}|_x$ is the fiber of ${\cal F}$ in $x$.
By Nakayama lemma, a coherent sheaf
is generated by any collection $\{s_i\}$ of sections which
generates ${\cal F}|_x={\cal F}\otimes_{\calo_M} (\calo_M/{\goth m}_x)$ for all $x\in M$.
Indeed, let ${\cal F}_x$ be the stalk of ${\cal F}$ in $x$.
Then ${\cal F}_x$ is a finitely-generated module
over a local ring $\calo_{M, x}$ of germs of holomorphic functions,
which is Noetherian by Lasker theorem (\cite[Chapter II, Theorem B.9]{_Gunning_Rossi_}).  
Nakayama lemma tells that a finitely-generated module $A$ over a Noetherian local
ring is generated by any set of elements which
generate $A$ modulo the maximal ideal.
Therefore, every stalk ${\cal F}_x$ of ${\cal F}$  is generated
by the images of $s_i$ in  ${\cal F}_x$; by definition,
this implies that ${\cal F}$ is generated by $\{s_i\}$.

We have shown that ${\cal F}$ is a vector bundle
whenever $\rk {\cal F}_x$ is constant in $x$.

Given a direct sum decomposition
${\cal F} = {\cal F}' \oplus {\cal F}''$,
we immediately obtain $\rk {\cal F}_x = \rk {\cal F}'_x + \rk {\cal F}''_x$.
If $\rk {\cal F}|_x$ is constant in $x$, this implies that 
$ \rk {\cal F}'_x =\const$ and $\rk {\cal F}''_x=\const$
because $\rk {\cal F}'_x$ is upper semicontinuous as
a function of $x$ by Nakayama lemma. This
implies that for any direct sum decomposition
${\cal F} = {\cal F}' \oplus {\cal F}''$
of a vector bundle onto a direct sum of coherent
sheaves, the summands are also vector bundles.
\endproof

\hfill

For the next proposition, we need the following claim,
which is almost trivial.

\hfill

\claim\label{_third_Mall_Claim_}
Let $0 \arrow A \arrow B \arrow C\arrow 0$
be an exact sequence of holomorphic bundles
on a Hopf manifold. Assume that $C$ and another of these bundles,
$B$ or $A$, are Mall. Then the third one is also
Mall.

\hfill

{\bf Proof:}
The functor $\pi^*$ is exact, and $j_*$ is left exact.
Therefore, applying the functor $j_*\pi^*$ to this sequence,
we obtain the following long exact sequence
\begin{equation}\label{_exact_exte_Equation_}
0 \arrow j_*\pi^*A \stackrel \alpha \arrow j_*\pi^*B \stackrel \beta \arrow j_*\pi^*C\arrow 
R^1 j_* \pi^* A^* \arrow ...
\end{equation}
where $R^ij_*$ is the higher derived pushforward functor.
If $C$ is Mall, then the sheaf  $j_*\pi^*C$ is locally free. Therefore,
the map $\beta$ has a section $\beta':\;j_*\pi^*C\arrow j_*\pi^*B $. This gives 
a subsheaf $\alpha(j_*\pi^*A)\oplus  \beta'(j_*\pi^*C) \subset  j_*\pi^*B$,
which by construction  coincides with
$j_*\pi^*B$ outside of 0. 
By Hartogs, it has to coincide with $j_*\pi^*B$ in 0 as well.
Then we have a direct sum decomposition
$j_*\pi^*A\oplus  \beta'(j_*\pi^*C) = j_*\pi^*B$.
Therefore, one of the sheaves $j_*\pi^*A$ and $j_*\pi^*B$
is locally free when the other one is locally free
(the sheaf $j_*\pi^*C=\beta'(j_*\pi^*C)$ is locally free by assumption).
\endproof

\hfill

The next Proposition contains a list of examples of Mall bundles.

\hfill

\proposition\label{_Mall_exa_Proposition_}
Let $H$ be a Hopf manifold. 
Then any line bundle on $H$, any tensor bundle 
$TH^{\otimes p}\otimes_{\calo_H} T^*H^{\otimes q}$,
their tensor products and direct sum components are Mall.

\hfill

\proof
Tensor products and direct sum components were treated 
in \ref{_tensor_direct_component_Mall_Claim_}.
Let $L$ be a line bundle on a Hopf manifold,
and $j_*\pi^*L$ the corresponding sheaf on $\C^n$
(\ref{_Mall_bundle_Definition_}). 
By Siu's extension theorem, \cite[Main Theorem]{_Siu:Extension_}, $j_*\pi^*L$ is 
a normal (and, therefore, reflexive) coherent sheaf of
rank 1, hence it is locally free
(\cite[Ch. II, Lemma 1.1.15]{_oss_}).
This implies that $L$ is Mall. 

The tangent bundle $TH$ is Mall because
$\pi^* TH$ is the tangent bundle $T(\C^n\backslash 0)$;
since the tangent bundle $T\C^n$ is trivial, 
$j_* \pi^* TH= T\C^n$ by Hartogs theorem.
\endproof

\subsection{$G$-equivariant sheaves}

\definition
Let $M$ be a topological space equipped
with the action of the group $G$ by continuous maps.
{\bf A $G$-equivariant sheaf} is a sheaf ${\cal F}$
equipped with a collection of isomorphisms
$R_g:\; {\cal F} \tilde \arrow g^*({\cal F})$
for all $g \in G$, which defines a $G$-action on the
\'etale space of ${\cal F}$.\footnote{This notion
is called {\bf a $G$-sheaf} in
\cite{_Grothendieck:Tohoku_}.}

\hfill

\remark Clearly, the equivariance of the $G$-action on 
${\cal F}$ can be translated into the {\bf equivariance 
relation} $R_{g_1g_2}= g_2^*(R_{g_1}) R_{g_2}$, for all
$g_1, g_2 \in G$.

\hfill

\theorem\label{_covering_equivariant_Theorem_}
Let $M$ be a locally connected, locally simply connected
topological space, and $\pi:\; M_1 \arrow M$ a Galois cover,
that is, a covering equipped with a free action of the
group $G$ such that $M = M_1/G$. 
Then the category of sheaves on $M$ is equivalent
to the category of $G$-equivariant sheaves\index{cover!Galois}
on $M_1$.

\hfill

\proof
Let ${\cal F}$ be a sheaf on $M$. Then the sheaf
$\pi^*({\cal F})$ is $G$-equivariant; indeed, 
$g\circ\pi= \pi$ for any $g\in G$, hence
$\pi^*({\cal F})= g^*(\pi^*({\cal F}))$,
which defines a $G$-equivariant structure on $\pi^*({\cal F})$.

Conversely, let ${\cal F}_1$ be a $G$-equivariant
sheaf on $M_1$,  $U\subset M$, and $U_1:=\pi^{-1}(U)$.
The equivariant structure on ${\cal F}_1$ defines
a $G$-action on the space of sections
${\cal F}_1(U_1)$. Indeed,
$g^*({\cal F}_1)(U_1)$ is by definition isomorphic to ${\cal F}_1(U_1)$,
and the isomorphisms $R_g:\; g^*({\cal F}_1)(U_1)\arrow {\cal F}_1(U_1)$
can be interpreted as self-maps on this space indexed by $g\in G$.
The equivariance relation $R_{g_1g_2}= g_2^*(R_{g_1}) R_{g_2}$
implies that the composition of these self-maps 
is compatible with the multiplication in $G$.

Consider the sheaf ${\cal F}_1^G$ on $M$ with the
space of sections ${\cal F}_1^G(U)$ equal to the space
of $G$-invariant sections of ${\cal F}_1(\pi^{-1}(U))$.
We claim that this functor from the category
of $G$-equivariant sheaves on $M_1$ to the
category of sheaves on $M$ is inverse to the
pullback functor defined above.

Clearly, any section of ${\cal F}$ on $U$ can be interpreted as
a $G$-invariant section of $\pi^*({\cal F})$ on $\pi^{-1}(U)$,
and any $G$-invariant section of $\pi^*({\cal F})$ comes from ${\cal F}(U)$.
Therefore, the composition of $\pi^*$ and $(\cdot )^G$
is equivalent to identity. To prove that
$\pi^*$ and $(\cdot )^G$ are inverse,
it remains to show that the functor
${\cal F}_1 \arrow \pi^*({\cal F}_1^G)$
is naturally isomorphic to identity 
(\cite[Section IV.4]{_Mac_Lane:Categories_}).
We leave this as an exercise to the reader.
\endproof

%
%
%


\section{Resonance in Mall bundles}


\subsection{Resonant matrices}

\definition\label{_resonant_matrix_Definition_}
Let $A\in \GL(n, \C)$ be a matrix with eigenvalues $\alpha_1,..., \alpha_n$.
This matrix is called {\bf resonant} if there exists
a relation $\alpha= \prod_{i=1}^n \alpha_i^{k_i}$,
with $\alpha$ an eigenvalue of $A$ and 
$k_i \in \Z^{\geq 0}$, $\sum_i k_i \geq 2$, and {\bf non-resonant}
otherwise. 

\hfill

The reason we need this definition is the following elementary lemma.

\hfill

\lemma\label{_inva_tensors_vanish_no_resonant_Lemma_}
Let $\gamma:\; \C^n \arrow \C^n$ be a germ of a holomorphic diffeomorphism in 0.
Let $B:= T\C^n\otimes_{\calo_{\C^n}} (T^* \C^n)^{\otimes k}$, $k >1$.
Diffeomorphisms of $\C^n$ can be naturally extended to $T\C^n$ and to its
tensor powers. Let $B_0$ denote the space of germs of sections of $B$ in 0;
clearly, $\gamma$ induces a natural automorphism of $B_0$. 
Assume that the differential $A:=D_0\gamma$ of $\gamma$ in 0 is non-resonant. 
Then any $\gamma$-invariant germ $v\in B_0$ vanishes.

\hfill

\proof
Let $t_1, ..., t_n$ be coordinates on $\C^n$.
Consider the Taylor series decomposition for 
$v$:
\[
v = \sum_i v_i  P_i(t_1, ..., t_n)
\]
where $P_i$ is a homogeneous polynomial of degree $i$
and $v_i \in V \otimes (V^*)^{\otimes k}$, where $V= T_0 \C^n$.
Let $l$ be the smallest integer such that $P_l \neq 0$.
The corresponding Taylor term can be considered as
an element of $V \otimes (V^*)^{\otimes k}\otimes \Sym^l(V^*)$.
By the chain rule, this tensor is coordinate-independent, hence
it is also $A$-invariant. Let $\alpha_1, ..., \alpha_n$
be the eigenvalues of $A$ on $T_0 \C^n$. Then 
$\alpha_1^{-1}, ..., \alpha_n^{-1}$ are the 
eigenvalues of $A$ on $T_0^*\C^n$. This is clear,
because $A$ preserves the pairing between $T\C^n$ and 
$T^*\C^n$: the differential $D\gamma$ acts covariantly
on vector field and contravariantly on 1-forms.

We obtain that the eigenvalues
of $A$ on $V \otimes (V^*)^{\otimes k}\otimes \Sym^l(V^*)$ are products
of the form $\alpha_l\prod_{i=1}^k \alpha_{j_i}^{-1}\prod_{i=1}^l \alpha_{s_i}^{-1}$.
Unless $A$ acts on $V \otimes (V^*)^{\otimes k}\otimes \Sym^l(V^*)$ 
with an eigenvalue 1, there would be no $A$-invariant vectors. Therefore
$v=0$ unless $\alpha_l = \prod_{i=1}^n \alpha_i^{k_i}$, where
$\sum_i k_i \geq 2$, for some eigenvalue $\alpha_l$.
\endproof

\subsection{Resonant equivariant bundles}

\definition
Let $\gamma$ be an invertible holomorphic contraction on $\C^n$,
centered in 0, and $B$ a $\gamma$-equivariant holomorphic
vector bundle on $\C^n$. Let $\alpha_1, ..., \alpha_n$
be the eigenvalues of the differential $A=D_0\gamma$ of
$\gamma$ in 0, and $\beta_1, ..., \beta_m$ the 
eigenvalues of the differential of the equivariant action of $\gamma$ 
on the fiber $B\restrict 0$. We say that $B$ is {\bf non-resonant}
if there is no multiplicative relation of the form $\beta_i = \beta_j \prod_{l=1}^k\alpha_{i_l}$
for some integer $k \geq 1$, where $\{\alpha_{i_1}, ...,\alpha_{i_k}\}$ is any collection of 
$k$ eigenvalues, possibly repeating, and $\beta_i, \beta_j\in \{\beta_1, ..., \beta_m\}$
any two eigenvalues, possibly the same. 

\hfill

\remark
A linear map with eigenvalues $\alpha_1, ..., \alpha_n$ is
resonant if one of the eigenvalues is a product of two or more
eigenvalues. The data associated with a $\gamma$-equivariant
bundle consists of linear operators, the differential $D\gamma\restrict {T_0\C^n}$
with the eigenvalues $\alpha_i$
and the differential  of the equivariant action of $\gamma$ 
on the fiber $B\restrict 0$, with the eigenvalues $\beta_j$.
It has resonance when one of $\beta_j$ is a product of another
and one or more $\alpha_i$'s. The resonance in the bundle
$T\C^n$ with the natural $\gamma$-equivariant structure
is the same as the resonance in the linear operator
$D\gamma\restrict {T_0\C^n}$.

In other words,
if $B$ is $T\C^n$ with the standard $\gamma$-equivariant structure,
the relation $\beta_i = \beta_j \prod_{l=1}^k\alpha_{i_l}$
becomes $\alpha_i = \alpha_j \prod_{l=1}^k\alpha_{i_l}$, $k \geq 1$,
that is, $B=T\C^n$ is non-resonant if and only if
the differential $A=D_0\gamma$ is non-resonant.

\hfill

Resonant automorphisms can be characterized in terms 
of invariant 1-forms with coefficients in 
endomorphisms of $B$, as follows.

\hfill

\theorem\label{_resonance_in_equiv_bundles_Theorem_}
Let $\gamma:\; \C^n \arrow \C^n$ be an invertible holomorphic
contraction centered in 0, and $B$ a $\gamma$-equivariant vector bundle.
Then $B$ is resonant if and only if there exists a non-zero
$\gamma$-invariant section of the bundle $\Omega^1\C^n\otimes \End(B)$.

\hfill

\pstep
Let $R$ be a $\gamma$-invariant section of $\Omega^1\C^n\otimes \End(B)$. 
We are going to prove that $B$ is resonant.
By \cite[Theorem 5.3.1]{_Forstneric:Oka_book_}, $B$ is trivial.
Choose a basis $b_1,..., b_m$ of $B$. Let $b_{ij}\in \End(B)$ be the
corresponding elementary matrices, defining a basis in $\End(B)$.
We write $R= \sum_{i,j,l} f_{ijl} dz_l \otimes b_{ij}$, where $z_i$ are coordinates in $\C^n$, and 
$f_{ijl}\in \calo_{\C^n}$ a function.
We write each $f_{ijl}$ as Taylor series, $f_{ijl} = \sum_s P_s^{ijl}$, 
where $P_s^{ijl}(z_1, ..., z_n)$ is 
a homogeneous polynomial of coordinate functions $z_i$
of degree $s$. Let $d$ be the smallest number
for which not all $P_d^{ijl}$ vanish. Using the chain rule
again, we obtain that 
$\sum_{i,j} P_d ^{ijl} dz_l \otimes b_{ij}$ is $D_0\gamma$-invariant. 
This is a $D_0\gamma$-invariant vector in the space
$W^* \otimes \End(B|_0) \otimes \Sym^d(W^*)$, where $W=T \C^n$, with the
action of $D_0\gamma$ which comes from the tensor product.

The eigenvalues of $D_0\gamma$ on the space
$W^* \otimes \End(B|_0) \otimes \Sym^d(W^*)$ are eigenvalues of
$D_0\gamma$ on $W^*$ times the eigenvalues on $\End(B|_0)$
times a product of $d$ eigenvalues of
$D_0\gamma$ on $W^*$, that is, a numbers of form 
$\beta_u \beta_v^{-1} \alpha_l^{-1}$. 

Since
$\sum_{i,j} P_d ^{ijl} dz_l \otimes b_{ij}\in W^* \otimes \End(B|_0) \otimes \Sym^d(W^*)$ 
is $D_0\gamma$-invariant, one of these numbers is  1. 
This gives a relation
$\beta_u  = \beta_v\alpha_l\prod_{l=1}^d\alpha_{i_l}$.

\hfill

{\bf Step 2:} Let $B$ be a resonant bundle;
we are going to produce an invariant section of 
$\Omega^1\C^n\otimes \End(B)$. 
We follow the standard scheme
(see \ref{_contra_compact_Theorem_},
\cite[Theorem 1.3]{ov_non_linear}),
using the Riesz-Schauder theorem 
and compactness of the action of holomorphic contractions on
holomorphic functions; this time, the contraction
acts on the sections of an equivariant bundle.
Let $M_0\subset \C^n$
be a subset which satisfies $\gamma(M_0) \Subset M_0$.
We equip $B$ with a Hermitian metric, and notice that
the space $H^0_b(M_0, \Omega^1M_0\otimes \End(B))$ with sup-norm
is Banach, by Montel theorem. Denote by 
$V$ the bundle $\Omega^1M_0\otimes \End(B)$.
By the standard argument (\ref{_contra_compact_Theorem_}), the operator
$\gamma^*:\; H^0_b(M_0, V)\arrow H^0_b(M_0, V)$ is compact.

Consider the filtration on $H^0_b(M_0, V)$
by the powers of the maximal ideal ${\goth m}_0$ of zero,
\[
H^0_b(M_0, V) \supset
H^0_b(M_0, {\goth m}_0V)\supset 
H^0_b(M_0, {\goth m}_0^2V)\supset ...
\]
The finite-dimensional space 
$\frac{H^0_b(M_0, V)}{H^0_b(M_0,  {\goth m}^k_0V)}$ is interpreted as
the space of $(k-1)$-jets of the sections of $V$ in 0.

Using the integral Cauchy formula,  any 
derivative of a function in a point 
can be expressed through a certain integral
of this function.
Therefore, the restriction map 
\[ H^0_b(M_0, V)\arrow
\frac{H^0_b(M_0, V)}{H^0_b(M_0,
{\goth m}^k_0V)}=  
V\restrict 0
\]
taking a section of $V$ to its $(k-1)$-jet is also
continuous in sup-norm. This implies also 
that the subspaces $H^0_b(M_0, {\goth m}_0^kV)\subset H^0_b(M_0, V)$
of sections with vanishing $(k-1)$-jet are closed.

The map $\gamma^*:\; H^0_b(M_0, V)\arrow H^0_b(M_0, V)$
preserves this filtration. Therefore, the eigenvalues of 
$\gamma^*$ on $H^0_b(M_0, V)$ and on the associated graded
space $\bigoplus_i \frac{H^0_b(M_0, {\goth m}^{i-1}V)}{H^0_b(M_0,
{\goth m}^i_0V)}$ are equal.
However, the space 
$\frac{H^0_b(M_0, {\goth m}^{d-1}V)}{H^0_b(M_0,
{\goth m}^d_0V)}$ is naturally identified with
$W^* \otimes_\C \Sym^d W^* \otimes_\C \End(B|_0)$,
where $W= T_0 \C^n$ (Step 1). In Step 1, we identified
this space with the space of $(d-1)$-th Taylor coefficients of
a section of $V$.

The eigenvalues of
$\gamma^*$ on this space are 
$\beta_1 \beta_2^{-1}\prod_{i=1}^{d+1}\alpha_i^{-1} $,
where $\beta_1, \beta_2$  denotes some eigenvalues of
the equivariant action of $\gamma^*$ on $B|_0$, possibly equal,  and
$\alpha_1, ..., \alpha_{d+1}$ is a collection of some
eigenvalues of $D_0\gamma$ on $W$, also possibly equal.
The existence of resonance on $B$ implies that
$\prod_{i=1}^{d+1} \alpha_i^{-1} \beta_1 \beta_2^{-1}=1$
for an appropriate choice of the eigenvalues
$\alpha_1, ..., \alpha_{d+1}, \beta_1, \beta_2$,
and this is equivalent to $\gamma^*$ having eigenvalue 1
on some $\gamma^*$-invariant quotient of $H^0_B(M_0, V)$, 
hence on $H^0_B(M_0, V)$, too.
\endproof

\hfill

We can rewrite this theorem as a result about Mall bundles.

\hfill

\corollary\label{_B_on_Hopf_resonant_via_sections_Corollary_}
Let $B_H$ be a Mall bundle on a Hopf manifold
$H= \frac{\C^n \backslash  0}{\langle \gamma \rangle}$,
and $B:=j_* \pi^* B_H$ the corresponding
$\gamma$-equivariant bundle on $\C^n$.
Then $B$ is resonant if and only if
$H^0(H, \Omega^1 H\otimes_{\calo_H}(\End B_H))\neq 0$.
\endproof

\hfill

When $B= \Omega^1 \C^n$, the equivariant bundle $B$
is resonant if and only if the action of $D_0\gamma^*$
on $T_0 \C^n$ is resonant:

\hfill

\corollary\label{_resonant_matrix_tangent_bundle_Corollary_}
Let $\gamma:\; \C^n \arrow \C^n$ be an invertible
contraction centered in zero, $W:= T_0\C^n$, and $A=D_0\gamma\in \End W$ 
its differential in zero. Then the following are equivalent.
\begin{description}
\item[(i)] The matrix $A$ is resonant.
\item[(ii)] There exists a non-zero $\gamma^*$-invariant section of
$\Omega^1(\C^n) \otimes \End(T\C^n)$.
\item[(iii)] The bundle $T\C^n$ with the natural 
equivariant structure induced by the action of $\gamma$ is
resonant.
\end{description}

\proof
From the definition it is clear that $A$ is resonant if
and only if the $\gamma$-equivariant bundle $B=\Omega^1 \C^n$
is resonant. By \ref{_resonance_in_equiv_bundles_Theorem_},
this is equivalent to the existence of non-zero $\gamma^*$-invariant sections of
$\Omega^1(\C^n) \otimes \End(T\C^n)$.
\endproof

\hfill

Consider an invertible
contraction $\gamma:\; \C^n \arrow \C^n$ centered in zero,
and let $B$ be a $\gamma$-equivariant vector bundle.
Let $R= (\Omega^1 \C^n)^{\otimes k}$ be the bundle of
$k$-multilinear forms on $\C^n$. We consider
$R$ as a $\gamma$-equivariant vector bundle as well.
Then the set of eigenvalues of the $\gamma^*$-action on
$H^0_b(M_0, R \otimes \End B)$ is equal to 
\[ 
\left\{\beta_1\beta_2^{-1}\prod_{i=1}^{d+k} \alpha_i^{-1} \right\}
\]
Here $d\in \Z^{\geq 0}$, $\alpha_1,..., \alpha_{d+k}$
is any collection of eigenvalues 
of $D_0\gamma^*$ on $W=T_0\C^n$, and $\beta_1, \beta_2$
some eigenvalues of $\gamma^*$ on $B|_0$.
This is proven by the same argument as proves
\ref{_resonance_in_equiv_bundles_Theorem_}, Step 2.
We obtained the following corollary.

\hfill

\corollary\label{_sections_tensor_w_coeff_in_End_B_Corollary_}
Let $\gamma:\; \C^n \arrow \C^n$ be an invertible
contraction of $\C^n$ centered in zero,
and $B$ a $\gamma$-equivariant vector bundle.
Assume that $B$ is non-resonant. Consider
$R= (\Omega^1 \C^n)^{\otimes k}$, $k \geq 1$,
as a $\gamma$-equivariant bundle.
Then the space of $\gamma^*$-invariant sections
of $R\otimes_{\calo_{\C^n}} B$ is empty.
\endproof


\subsection{Holomorphic connections on vector bundles}


To go on, we need the notion of 
a holomorphic connection.

\hfill

\definition\label{_holo_conne_Definition_}
Let $B$ be a holomorphic vector bundle on 
a complex manifold, and $\nabla:\; B \arrow B\otimes_{\calo_M} \Omega^1 M$
a differential operator which satisfies $\nabla(fb) = df \otimes b + f \nabla B$
for any locally defined holomorphic function $f$ and any local section $b$ of $B$.
Then $\nabla$ is called {\bf a holomorphic connection}.

\hfill

The notion of a holomorphic connection was introduced by
M. Atiyah in \cite{_Atiyah:hol_bundles_}; for more results and
references, see
\cite{_Biswas_,_Biswas_Dumitrescu:holomorphic_affine_,%
_Dumitrescu_Biswas:holomorphic_Riemannian_}. 

\hfill

Every flat connection is holomorphic with respect
to the holomorphic structure induced by this connection,
but there are more holomorphic connections than there
are flat connections.  Indeed, holomorphic connections
can be realized as objects of differential geometry,
as follows.

\hfill

Let $(B, \nabla)$ be a complex vector bundle with connection on a 
complex manifold, and $\bar\6= \nabla^{0,1}$ the
corresponding $\bar\6$-operator.
By Koszul-Malgrange theorem (\cite[Chapter I, Proposition
3.7.]{_Kobayashi_Bundles_},
\cite{_Koszul_Malgrange_}),
$\bar\6$ defines a holomorphic structure on $B$ if and only
if $\bar\6^2=0$, or, equivalently, when 
$\nabla^2 \in [\Lambda^{2,0}(M)\oplus \Lambda^{1,1}(M)]\otimes \End(B)$.

\hfill

\proposition\label{_holo_conne_curvature_Proposition_}
Let $(B, \nabla)$ be a complex vector bundle with connection on a 
complex manifold, and $\bar\6= \nabla^{0,1}$ the corresponding $\bar\6$-operator.
Assume that $\bar\6^2=0$, and let ${\cal B}=\ker \bar\6$ be 
the holomorphic vector bundle obtained from $\bar\6$ 
using \cite[Chapter I, Proposition
3.7.]{_Kobayashi_Bundles_},
\cite{_Koszul_Malgrange_}. Then the
following assertions are equivalent.
\begin{description}
\item[(i)] The operator $\nabla^{1,0}$ is a holomorphic
connection operator on ${\cal B}$.
\item[(ii)] $\nabla^2 \in \Lambda^{2,0}(M)\otimes \End(B)$.
\end{description}
\proof
If $\nabla^{1,0}$ is a holomorphic
connection operator, it commutes with $\bar\6$, hence
the $(1,1)$-part of the curvature of $\nabla$ vanishes.
Conversely, if $\nabla^2 \in \Lambda^{2,0}(M)\otimes \End(B)$,
this implies that $\nabla^{1,0}$ commutes with $\bar\6$,
hence $\nabla^{1,0}$ maps holomorphic sections of $B$
to the holomorphic sections of $B \otimes \Lambda^{1,0}(M)$.
This implies that $\nabla^{1,0}$ is a holomorphic connection.
\endproof

\hfill

\remark
A holomorphic vector bundle $B$
has Chern classes represented by closed forms
of type $(p,p)$, because it admits Chern connections. However, 
the Chern classes of a bundle equipped with a holomorphic
connection are represented by holomorphic forms of type $(2p,0)$.
On a K\"ahler manifold (or any other manifold admitting
the $(p,q)$-decomposition in cohomology), this is impossible,
unless $c_i(B)=0$ for all $i$. This is why holomorphic
connections rarely occur in K\"ahler geometry. However,
on non-K\"ahler manifolds they don't seem that rare.

\hfill

The Picard group $\Pic(M)$ of a complex manifold  $M$ (that is, the group
of line bundles, with the group operation defined by the
tensor multiplication) is naturally identified with  
$H^1(M, \calo_M^*)$. Similarly, one could indentify
the group of 1-dimensional local systems with $H^1(M, \C_M^*)$.
Here $\C^*_M$
denotes the constant sheaf with the space of sections
$H^0(U, \C^*_M)=\C^*$ for each connected
open set $U\subset M$.  The natural map
$H^1(M, \C_M^*)\arrow H^1(M, \calo_M^*)$
can be interpreted as a forgetful map,
taking a flat line bundle $(L, \nabla)$
to the holomorphic line bundle $(L, \nabla^{0,1})$.

\hfill

For K\"ahler manifolds, the map 
$H^1(M, \C_M^*)\arrow H^1(M, \calo_M^*)$ is never an isomorphism,
unless $b_1(M)=0$. Indeed, from the exponential exact sequence
it follows that $\dim_\C H^1(M, \C_M^*)= b_1(M)$ and
$\dim_\C H^1(M, \calo_M^*)= \frac 1 2 b_1(M)$.
However, on a Hopf manifold, these two groups are isomorphic,
that is, every holomorphic line bundle admits a unique flat connection.

\hfill

For the sake of completeness, 
we give the proof of the following result of Kodaira:

\hfill

\proposition\label{_Picard_Hopf_Proposition_}
(\cite[page 57]{_Kodaira_Structure_III_}) \\ Let $H$ be a Hopf manifold,
and $H^1(H, \C_H^*)$ the cohomology with coefficients in
the constant sheaf $\C_H^*$.
Then the natural map $H^1(H, \C_H^*)\arrow H^1(H, \calo_H^*)$
to the Picard group is an isomorphism.

\hfill

\pstep
Consider the exponential exact sequences
\[
0 \arrow \Z_H \arrow \calo_H \arrow \calo_H^* \arrow 0
\]
and 
\[
0 \arrow \Z_H \arrow \C_H \arrow \C_H^* \arrow 0.
\]
The corresponding long exact sequences of cohomology give
\[
0 \arrow H^1(\Z_H) \arrow H^1(\calo_H)\arrow H^1(\calo_H^*)\arrow H^2(\Z_H)=0
\]
and 
\[
0 \arrow H^1(\Z_H) \arrow H^1(\C_H)\arrow H^1(\C_H^*)\arrow H^2(\Z_H)=0.
\]
It remains to show that the natural map 
$\nu:\; H^1(\C_H)\arrow H^1(\calo_H)$
is an isomorphism; however, both groups are equal to $\C$
(\ref{_Dolbeault_for_Hopf_Theorem_}), hence it would
suffice to show that $\nu \neq 0$.
The relevant $E_2$-term of the Dolbeault spectral sequence
is $(H^{0,1}(H)=\C; H^{1,0}(H)=0)$ 
(\ref{_Dolbeault_for_Hopf_Theorem_}, \ref{_holo_forms_Hopf_vanish_Proposition_}).
The sum of dimensions of these spaces is 1. However, $b_1(H)=1$,
hence the higher differentials of this spectral sequence  vanish
on $E_2^{1,0}+E_2^{0,1}$ and it degenerates in the $E_2^{1,0}+E_2^{0,1}$-term.
This implies that the first de Rham cohomology $H^1(H)$
is equal to the first Dolbeault cohomology
$H^{0,1}(H)\oplus H^{1,0}(H)$.

\hfill

{\bf Step 2:}
The standard map
$\nu:\; H^1(\C_H)\arrow H^1(\calo_H)$ is
the natural map from the $E_\infty^{0,1}$-term of this
spectral sequence to $E_2^{0,1}$, which is an
isomorphism because the spectral sequence degenerates.
\endproof


\subsection{The flat connection on a non-resonant Mall bundle}


\definition
Let $B_H$ be a Mall bundle on a Hopf manifold
$H= \frac{\C^n \backslash  0}{\langle \gamma \rangle}$,
and $B:=j_* \pi^* B_H$ the corresponding
$\gamma$-equivariant bundle on $\C^n$.
We call $B_H$ {\bf resonant}
if $H^0(H, \Omega^1_H\otimes \End B_H)\neq 0$,
or, equivalently, when $B=j_* \pi^* B_H$
is a resonant $\gamma$-equivariant bundle on $\C^n$
(\ref{_B_on_Hopf_resonant_via_sections_Corollary_}).

\hfill

\remark\label{_conne_unique_Remark_}
Let $B_H$ be a non-resonant Mall bundle on a Hopf manifold.
Since the difference of two connections is
a holomorphic 1-form with coefficients in endomorphisms,
and $H^0(H, \Omega^1_H\otimes \End B_H)=0$,
a holomorphic connection on $B_H$ is unique, if it exists.

\hfill

\theorem\label{_Mall_flat_connection_Theorem_}
Let $B$ be a holomorphic vector bundle over a Hopf manifold $H$,
$\dim_\C H \geq 3$.  
Assume that $B$ admits a flat connection compatible
with the holomorphic structure. Then $B$ is Mall.
Conversely, any non-resonant Mall bundle on $H$
admits a flat connection $\nabla$.\footnote{All flat connections
are holomorphic (\ref{_holo_conne_curvature_Proposition_}). 
By \ref{_conne_unique_Remark_}, $B$
admits a unique holomorphic connection, hence $\nabla$ is unique.}

\hfill

\pstep 
Let $(B, \nabla)$ be a holomorphic bundle with a flat connection on $H$,
and $\pi:\; \C^n \backslash 0 \arrow H$ the
universal cover. Since $\pi_1(\C^n \backslash 0)=0$,
the flat bundle $\pi^* B$ is trivial, hence $B$ is Mall.

\hfill

{\bf Step 2:} 
Let $B$ be a non-resonant Mall bundle on a Hopf manifold $H$.
We are going to show that $B$ admits a holomorphic connection.
Locally, a holomorphic connection always exists, and the
difference between two holomorphic connections is a
section of $\Omega^1 U
\otimes_{\calo_U} \End (B)$. Therefore, the 
obstruction to the  existence 
of a holomorphic connection belongs to 
$H^1(H, \Omega^1 H\otimes_{\calo_H} \End (B))$. 
Since $B$ is Mall, we have
\[ \rk H^1(H, \Omega^1 H\otimes_{\calo_H} \End (B))= 
\rk H^0(H, \Omega^1 H\otimes_{\calo_H} \End (B))
\]
(\ref{_Mall_cohomology_Theorem_}). 
When $B$ has no resonance, 
$\rk H^0(H, \Omega^1 H\otimes_{\calo_H} \End (B))=0$
hence $B$ admits a connection.

\hfill

{\bf Step 3:} To finish
\ref{_Mall_flat_connection_Theorem_},
it remains to show that any holomorphic connection
$\nabla$ on a non-resonant Mall bundle $B$ is flat. 
However, its curvature is a holomorphic 2-form, 
and the space of holomorphic 2-forms with
coefficients in $\End(B)$ vanishes by 
\ref{_sections_tensor_w_coeff_in_End_B_Corollary_}.
\endproof


\section{Flat connections on Hopf manifolds}
\label{_flat_Hopf_Section_}

\subsection{Developing map for flat affine manifolds}

For an introduction to flat affine manifolds, see 
\cite{_Goldman:structures_book_,_Abels:survey_,_Shima:book_} and the references therein.
Recall that {\bf an affine function} on a vector space is a linear function
plus constant.

\hfill

\definition
Let $M$ be a manifold equipped with a sheaf ${\cal F}\subset \C^\infty M$.
We say that ${\cal F}$ {\bf defines a flat affine structure on $M$}
if for each $x\in M$ there exists a neighbourhood diffeomorphic
to the ball $B \subset \R^n$ such that ${\cal F}$ corresponds
to affine functions on $B$. In other words, $M$ is a flat
affine manifold if there exists an open cover $\{ U_i\}$
with all $U_i$ diffeomorphic to an open ball in $\R^n$ and
all transition maps are affine. 
Such a cover is called {\bf an affine atlas} of $M$.
The sheaf ${\cal F}$ is called {\bf the sheaf of affine functions}
on the flat affine manifold $M$.

\hfill

Flat affine structures can be equivalently described in terms
of torsion-free, flat connections. The following
proposition is well-known; we include
its proof for  completeness.

\hfill

\proposition\label{_flat_affine_via_connection_Proposition_}
Let $(M, {\cal F})$ be a flat affine manifold. Then
$M$ admits a unique torsion-free, flat connection $\nabla$ such that
the sections $f$ of ${\cal F}$ satisfy $\nabla( df)=0$. Conversely,
if $\nabla$ is a torsion-free flat connection on $M$,
then the sheaf $\{ f \in C^\infty M\ \ |\ \ \nabla(df)=0\}$
defines a flat affine structure.

\hfill

\proof
Let $(M, {\cal F})$ be a flat affine manifold, and
$\{ U_i\}$ its affine atlas. Each $U_i$ admits 
a connection $\nabla:\; TM \arrow TM \otimes \Lambda^1 M$ with 
\begin{equation}\label{_standard_flat_connection_Equation_} 
\nabla\left(\sum_j f_j \frac{\6}{\6x_j}\right) = \sum_j df_j \otimes \frac{\6}{\6x_j}
\end{equation}
where $x_j$ are coordinate functions; this connection is clearly
flat and torsion-free, and $\ker (\nabla d)$ is the sheaf of affine
functions on $U_i$. We call \eqref{_standard_flat_connection_Equation_}
{\bf the standard flat connection} on $U\subset \R^n$.

Conversely, let $\nabla$ be a torsion-free affine connection on $M$.
Locally in an open subset $U\subset M$, 
the bundle $\Lambda^1M$ admits a basis $\theta_1, ..., \theta_n$ 
of $\nabla$-parallel sections,
which are closed because $\nabla$ is torsion-free (here we use
the relation $d\theta = \Alt(\nabla \theta)$, which holds
for any torsion-free connection).  This implies that
$\theta_i = dx_i$ whenever $U$ is simply connected.
The functions $x_i$ give a coordinate system on $U$,
because $dx_i$ are linearly independent, and in this
coordinate system the affine functions are those
which satisfy $\nabla (df)=0$.
\endproof

\hfill

Further on, we will also call a torsion-free flat connection
{\bf an affine structure} and a pair $(M, \nabla)$ 
{\bf an affine manifold}, or {\bf a flat affine manifold}.

\hfill

\definition
Let $\nabla$ be a connection on $TM$, where $M$ is a smooth manifold. 
We say that $(M, \nabla)$ is {\bf complete} if for any
$x\in M$ and any $v\in T_x M$, there exists a solution
$\gamma_v:\; \R \arrow M$ of the geodesic equation 
$\nabla_{\gamma_v'(t)}\gamma_v'(t)=0$ for all $t\in ]-\infty, \infty[$. 
{\bf The exponential map} is the map $\exp_x:\; T_x M \arrow M$
taking $v \in T_x M$ to $\gamma_v(1)$. Clearly, $\exp_x$ is a
diffeomorphism in a neighbourhood $U$ of $0\in T_x M$.
{\bf The developing map} $\dev:\; \exp_x(U) \arrow U$
is the inverse of $\exp_x$; in general, it is
defined only in a neighbourhood of $x\in M$.

\hfill

The following classical theorem, due to Auslander and
Markus, is well known. It can be considered as an
alternative definition of the complete flat affine
manifold. 

\hfill

\theorem\label{_dev_for_complete_affine_Theorem_} 
(\cite{_Auslander_Markus:holonomy_}, \cite{_Goldman:structures_book_}) 
Let $(M, \nabla)$ be a simply connected, connected  flat affine manifold.
Then the developing map can be extended to an affine map
$\dev:\; M \arrow T_x M$, also called {\bf the developing map}.
If, in addition, $M$ is complete, the developing map is 
a diffeomorphism. \endproof

%
%
%

\subsection{Flat affine connections on a Hopf manifold}

The following result will be used further on in 
our proof of the classical Poincar\'e linearization
theorem.

\hfill

\theorem\label{_flat_Hopf_is_linear_Theorem_}
Let $H= \C^n \backslash 0/\Z$ be a Hopf manifold.
Assume that $TH$ admits a torsion-free, flat connection
$\nabla$ compatible with the holomorphic structure. Then
$H$ is isomorphic to a linear Hopf manifold 
$\C^n \backslash 0/\langle A \rangle$, where $A\in \GL(n, \C)$ is a
linear endomorphism.

\hfill

\pstep
Let $\pi:\; \C^n \backslash 0\arrow H$ be the covering
map, and $j:\; \C^n \backslash 0\hookrightarrow \C^n$ 
the tautological open embedding map. Let $B$ be a tensor
bundle on $H$. By \ref{_Mall_exa_Proposition_}, $B$ is a Mall 
bundle, hence  $j_* \pi^* B$ is a bundle on $\C^n$.
By \cite[Theorem 5.3.1]{_Forstneric:Oka_book_}, 
any vector bundle on $\C^n$ is trivial.
Let $\nabla_0$ be the trivial connection on $T\C^n$.
Then the connection form $\nabla_0 -\pi^* \nabla\in \Omega^1 \C^n
\otimes_{\calo_{\C^n}} \End (T\C^n)$ is a holomorphic section of 
the trivial bundle $\Omega^1 \C^n\otimes_{\calo_{\C^n}} \End (T\C^n)$
over $\C^n \backslash 0$, hence it can be extended to a
holomorphic connection on $\C^n$. This connection
is also torsion-free, flat and holomorphic, by continuity
of the corresponding tensors. 

To avoid further confusion, we denote
$\C^n$ by $M$, and this connection by $\nabla_M$.
To finish the proof
of \ref{_flat_Hopf_is_linear_Theorem_}, it would
suffice to show that the flat affine manifold
$(M,\nabla_M)$ is isomorphic to $\C^n$ with the standard flat
connection. Since the $\Z$-action on $M$ preserves
the affine structure and fixes 0, it is defined by
a linear endomorphism, and $(M \backslash 0)/\Z$
is a linear Hopf manifold.

We conclude that it remains to show that the developing map
$\dev:\; M \arrow \C^n$ associated to $\nabla_M$ is an
isomorphism. This would follow from
\ref{_dev_for_complete_affine_Theorem_}
if we could prove that $\nabla_M$ is complete,
but we have no direct control over the connection
form of $\nabla_M$, hence the completeness is not obvious.

\hfill

{\bf Step 2:}
Let $\exp_0:\; U \arrow M$ 
be the exponential map, with $U\subset T_0M$
being a maximal neighbourhood of $0$ where $\exp_0$ is defined.

It remains to show that
$U=T_0M$; this would imply
that $(M, \nabla_M)$ is complete and the developing
map is an isomorphism. 

Denote by $A$ the generator
of the $\Z$-action on $M$ contracting $M$ to 0,
and let $A_0\in \End(T_0M)$ be its differential.
By construction, $\nabla_M$ is $A$-invariant. 
Therefore, $A$ preserves the geodesics of $(M, \nabla_M)$,
and we have a commutative diagram 
\[
\begin{CD}
U@>{\exp_0}>> M=\C^n \\
@V{A_0}VV  @VV{A} V \\
U @>{\exp_0}>> M= \C^n
\end{CD}
\]
This implies that $U\subset T_0 M$ is an $\Z$-invariant
neighbourhood of $0\in T_0 M$, where the action of $\Z=\langle A_0\rangle$ is
generated by $A_0$. Since 
$\bigcup_n A_0^{-n}(V)= T_0 M$ for any neighbourhood 
$V$ of 0, any $\langle A_0\rangle$-invariant neighbourhood
is equal to $T_0 M$. This implies that $(M, \nabla)$ is
complete, and the developing map is an isomorphism.
\endproof

\subsection{A new proof of Poincar\'e 
theorem about linearization of non-resonant contractions}

The Poincar\'e-Dulac theorem (\cite{_Arnold:ODE+_})
gives a normal form of a 
smooth (or analytic) contraction; its non-resonant case
is sometimes called {\em the Poincar\'e theorem}. It
proves that a contraction (or a germ of a contraction),
which is non-resonant, becomes linear after an appropriate
coordinate change. We give a new 
proof of this theorem based on complex geometry.
Note that the assumption $n \geq 3$ below is unnecessary;
we leave the case $n=2$ for the reader as an exercise.

\hfill

\theorem\label{_Poincare_new_proof_Theorem_}
Let $\gamma$ be an invertible holomorphic contraction of $\C^n$ centered in 0, 
$n \geq 3$. Assume that the differential
$D_0\gamma\in \GL(T_0 \C^n)$ is non-resonant. Then there exists
a holomorphic diffeomorphism $U:\; \C^n \arrow \C^n$ such that
$U \gamma U^{-1}$ is linear.\footnote{The same statement can
be stated for a germ of a holomorphic diffeomorphism; the proof will be essentially
the same. We leave the required changes to the reader.}

\hfill

\proof
Let $H:= \C^n \backslash 0 /\langle \gamma \rangle$ be the Hopf
manifold associated with $\gamma$, and 
$\pi:\; \C^n \backslash 0\arrow H$ the universal covering map.
By \ref{_Mall_exa_Proposition_}, the tangent
bundle $TH$ is Mall. By \ref{_resonant_matrix_tangent_bundle_Corollary_}, 
it is non-resonant.
By \ref{_Mall_flat_connection_Theorem_},
$TH$ admits a flat 
holomorphic connection $\nabla$. Since $\Omega^2 H \otimes TH$ is a direct sum
component of $\Omega^1 H \otimes \End(TH)$,
and $TH$ is non-resonant, 
\[ H^0(H,\Omega^2 H \otimes TH)
\subset H^0(H,\Omega^1 H \otimes \End(TH))=0,
\] 
(\ref{_resonant_matrix_tangent_bundle_Corollary_}),
hence $\nabla$ is torsion-free.

By \ref{_flat_Hopf_is_linear_Theorem_},
the universal cover $\C^n \backslash 0$ of $H$ admits flat coordinates
such that $\gamma$ is linear in these coordinates. 
This proves \ref{_Poincare_new_proof_Theorem_}.
\endproof

\hfill



\hfill

{\small

\noindent {\sc Liviu Ornea\\
University of Bucharest, Faculty of Mathematics and Informatics, \\14
Academiei str., 70109 Bucharest, Romania}, and:\\
{\sc Institute of Mathematics ``Simion Stoilow" of the Romanian
Academy,\\
21, Calea Grivitei Str.
010702-Bucharest, Romania\\
\tt lornea@fmi.unibuc.ro,   liviu.ornea@imar.ro}

\hfill

\noindent {\sc Misha Verbitsky\\
{\sc Instituto Nacional de Matem\'atica Pura e
Aplicada (IMPA) \\ Estrada Dona Castorina, 110\\
Jardim Bot\^anico, CEP 22460-320\\
Rio de Janeiro, RJ - Brasil }\\
also:\\
Laboratory of Algebraic Geometry, \\
Faculty of Mathematics, National Research University 
Higher School of Economics,
6 Usacheva Str. Moscow, Russia}\\
\tt verbit@verbit.ru, verbit@impa.br }


{\small 

\begin{thebibliography}{100}

\bibitem[Ab]{_Abels:survey_} H. Abels, {\em Properly discontinuous groups of affine transformations: a survey}, Geom. Dedicata {\bf 8}7 (2001), no. 1-3, 309-333.

\bibitem[Ar]{_Arnold:ODE+_} V. I. Arnol'd, {\em Geometrical
  Methods in the Theory of Ordinary Differential
  Equations,} Grundlehren der mathematischen Wissenschaften
  {\bf 250}, Springer, 1996.


\bibitem[At]{_Atiyah:hol_bundles_} M. F. Atiyah, {\em
  Complex analytic connections in fibre bundles}, Trans.
  Amer. Math. Soc. {\bf 85}(1), (1957) 181-207.

\bibitem[AM]{_Auslander_Markus:holonomy_} L. Auslander,
  L. Markus, {\em Holonomy of flat affinely connected
    manifolds}, Ann. of Math. (2) {\bf 62} (1955),
  139-151.


\bibitem[B]{_Biswas_} I. Biswas, {\em Vector bundles with
  holomorphic connection over a projective manifold with
  tangent bundle of nonnegative degree},
  Proc. Amer. Math. Math. {\bf 126} (1998), 2827-2834.

\bibitem[BD1]{_Biswas_Dumitrescu:holomorphic_affine_}
  I. Biswas, S. Dumitrescu, {\em Holomorphic affine
    connections on non-K\"ahler manifolds},
  Internat. J. Math. {\bf 27} (2016), no. 11, 1650094, 14
  pp.

\bibitem[BD2]{_Dumitrescu_Biswas:holomorphic_Riemannian_}
  I. Biswas, S. Dumitrescu, {\em Holomorphic Riemannian
    metric and the fundamental group},
  Bull. Soc. Math. France {\bf 147} (2019), no. 3,
  455-468.

\bibitem[D]{_Dulac_} H. Dulac, {\em Recherches sur les
  points singuliers des equations diff\'erentielles}
  J. Ecole Polytechn. Ser. II , {\bf 9} (1904),  1-25.


\bibitem[EnMath]{_Enc_Math:Poincare_Dulac_}
{\em Poincar\'e-Dulac theorem,} {Encyclopedia of Mathematics,}
\url{https://encyclopediaofmath.org/wiki/Poincar\%C3\%A9-Dulac_theorem}

\bibitem[Fo]{_Forstneric:Oka_book_}
F. Forstneri\v c, Stein Manifolds and Holomorphic Mappings, Springer Verlag, 2011.

\bibitem[Fr]{friedman} A. Friedman, Foundations of modern analysis, Dover, 2010.

\bibitem[GZ]{_Gan_Zhou_} 
N. Gan, X. Y. Zhou, {\em The
  cohomology of bundles on general non primary Hopf
  manifolds}. In: Recent Progress on Some Problems in
  Several Complex Variables and Partial Differential
  Equations. Contemp Math, Vol. {\bf 400}, 
Providence, RI: Amer Math Soc, 2006, 107-115.



\bibitem[G]{_Goldman:structures_book_} W. Goldman, 
Geometric Structures on Manifolds, AMS Open Math. Notes,
\url{https://www.ams.org/open-math-notes/omn-view-listing?listingId=111282}

\bibitem[Gro]{_Grothendieck:Tohoku_} 
A. Grothendieck,  {\em Sur quelques points 
d'alg\`ebre homologique, I, II}, Tohoku Math. J. (2) 
{\bf 9} (2, 3): 119-221. English translation: 
\url{http://www.math.mcgill.ca/barr/papers/gk.pdf}.


\bibitem[GN]{_Gunning_Rossi_} R. C. Gunning, H. Rossi, 
{\em Analytic functions of several complex variables,} 
Reprint of the 1965 original. AMS Chelsea Publishing, Providence, RI, 2009.

\bibitem[Ha]{_Haefliger:Hopf_} 
A. Haefliger, {\em Deformations of transversely 
holomorphic flows on spheres and deformations of Hopf manifolds}, 
	Compositio Math. {\bf 55} (1985), no. 2, 241-251. 


\bibitem[I]{_Ise_} M. Ise, {\em On the geometry of Hopf manifolds}, Osaka J. Math. {\bf 12} (1960), 387-402.

\bibitem[Kob]{_Kobayashi_Bundles_} S. Kobayashi, 
{\em Differential Geometry of Complex Vector Bundles,} 
Princeton Univ. Press, 1987.

\bibitem[Kod1]{_Kodaira_Structure_II_} K. Kodaira, {\em On the structure of compact complex surfaces, II}, Amer. J. Math. {\bf 88} (1966), 682-721.

\bibitem[Kod2]{_Kodaira_Structure_III_} K. Kodaira, {\em On the structure of compact complex surfaces, III}, Amer. J. Math. {\bf 90} (1968), 55-83.

\bibitem[KM]{_Koszul_Malgrange_} J.-L. Koszul, B. Malgrange, {\em Sur certaines
structures fibr\'ees complexes}, Arch. Math. {\bf 9} (1958), 102-109.

\bibitem[L]{_Lattes_} S. Latt\`es, {\em Sur les formes r\'eduites des transformations ponctuelles dans le domaine d'un point double}, Bull. Soc. Math. France {\bf 39} (1911), 309-345. 

\bibitem[Lib]{_Libgober_} A. Libgober, {\em Cohomology of
  bundles on homological Hopf manifolds},  Sci. China
  Ser. A {\bf 52} (2009), no. 12, 2688-2698. 

\bibitem[McL]{_Mac_Lane:Categories_} S. Mac Lane, 
{\em Categories for the working mathematician}, Springer, 1971.


\bibitem[M1]{_Mall_} D. Mall, {\em The cohomology of line bundles on Hopf manifolds}, Osaka J. Math. {\bf 28} (1991), 999-1015.

\bibitem[M2]{_Mall:Contractions_} D. Mall, {\em Contractions, Fredholm operators and the cohomology
of vector bundles on Hopf manifolds}, Arch. Math., {\bf 66} (1996), 71-76.

\bibitem[Mov]{_MO:Hopf_cohomology_}  
\url{https://mathoverflow.net/questions/25723/dolbeault-cohomology-of-hopf-manifolds}

\bibitem[OSS]{_oss_} C. Okonek, M. Schneider, H. Spindler,
  {\em Vector bundles on complex projective spaces}, 
Corrected reprint of the 1988 edition. With an appendix by
S. I. Gelfand. Modern Birkh\"auser
Classics. Birkh\"auser/Springer Basel AG, Basel, 2011. 

\bibitem[OV1]{ov_indam} L. Ornea, M. Verbitsky, {\em Embedding of LCK manifolds with potential into Hopf manifolds using Riesz-Schauder theorem}, ``Complex and Symplectic Geometry'', Springer INdAM series, 2017, 137-148.

\bibitem[OV2]{ov_non_linear} L. Ornea, M. Verbitsky, {\em
  Non linear Hopf manifolds are locally conformally
  K\"ahler}, arXiv:2202.12398.

\bibitem[RS]{_Ramani_Sankaran_} V. Ramani, P. Sankaran,
  {\em Dolbeault cohomology of compact complex homogeneous
    manifolds} Proc. Indian Acad. Sci. Math. Sci. {\bf
    109} (1999), no. 1, 11-21.

\bibitem[Sh]{_Shima:book_} H. Shima, {\em The geometry of
  Hessian structures}, World Scientific Publishing, 2007.


\bibitem[Si]{_Siu:Extension_} Y.-T. Siu, 
{\em Extension of locally free analytic sheaves}, Math. Ann. {\bf 179} (1969) 285-294.

\bibitem[St]{_Sternberg_contraction_} S. Sternberg, {\em
  Local contractions and a theorem of Poincar\'e},
  Amer. J.  Math. {\bf 79} (1957), 809-824.

\bibitem[Wu]{_Wu:Montel_} H. Wu, 
{\em Normal families of holomorphic mappings},
Acta Math. {\bf 119} (1967), 193-233.

\end{thebibliography}

}
\end{document}